# cSAT problem lower bound

*Abstract*—This article deals with the lower bound that is considered as the worst case minimal amount of time required to calculate a problem result for cSAT (counted Boolean satisfiability problem). It uses the observation that Boolean algebra is a complete first-order theory where every sentence is decidable. Lower bound of this decidability is defined and shown.

The article shows that deterministic calculation model made up of finite number of machines (algorithms), oracles, axioms, or predicates is incapable of solving considered NP-complete problem when its instance grows to infinity. This is a direct proof of the fact that P and NP complexity classes differ and oracle capable of solving NP-complete problems in polynomial time must consist of infinite number of objects (i.e., must be nondeterministic).

Corollary of this article clears complexity hierarchy:
**P < NP**

*Index terms*—complexity class, P vs NP, Boolean algebra, first order theory, first order predicate calculus.

## I. INTRODUCTION

Unknown relation between P and NP [5] complexity classes remains one of the significant unsolved problems in complexity theory. P complexity class consists of problems solvable by deterministic Turing machine (DTM) in polynomially bounded time, while NP complexity class consists of problem solvable by nondeterministic Turing machine (NDTM) in polynomially bounded time. This means that DTM can verify the solution of every NP problem in polynomially bounded time, even if polynomial algorithm for finding this solution is unknown [13].

All known attempts to prove whether these classes are or are not equal could not convince the community that arguments used there are final. Problem with attempts showing that P=NP is mainly with counter examples provided for methods described by solvers (see for example: [6], [9]), especially for large instances. Problem with proof attempts that P≠NP touches mainly the difference between a problem and an algorithm. Proving the inequality of these classes is equivalent to proving that "there is no such algorithm that solves a particular NP problem in polynomially bounded time." Algorithm is an immaterial object, so proving that it does not exist is rather difficult.

Can then the inequality of complexity classes be proved? One of the possible ways is to use the properties of first-order theory. Useful properties include every sentence φ in theory T is provable if there exists a set of axioms a, b, c, ... such that φ can be obtained using these axioms and the inference rules "modus ponens" and "universal generalization" (a ∧ b ∧ c... → φ) [2].

## II. BACKGROUND

This section presents some background for the first-order theory and other rules used in the article.

### A. First-Order Theories

First-order theory is a given set of axioms in some language. Language consists of logical symbols and set constants, functions, and relation symbols (predicates). Terms and formulas are built from language and give rise to sentences, which are formulas with no free variables in body.

Theory is then a set of sentences which may be closed if it contains all consequences of its elements. Theory can be also complete (i.e. every sentence can be proved or disproved), consistent (not every sentence is provable), or decidable (every sentence can be proved or disproved and there exists a computational path (algorithm) showing which sentences are provable).

An example of first-order theory that is complete and decidable is Boolean algebra [16] or Zermelo–Frænkel set theory.

### B. First-Order Logic

First-order logic, also called first-order predicate calculus (FOPC), is a system of deductions extending propositional logic. Atomic sentences of first-order logic are called predicates and are written usually in the form $P(t_1, t_2, ..., t_n)$. An important ingredient of the first-order logic not found in propositional logic is quantification.

In 1929, Gödel [8] proved that every valid logical formula is valid in first-order logic. In other words, it is proved that for complete first-order theory, inference rules of FOPC are sufficient to prove any valid formula.

First-order predicate calculus language consists of predicates, constants, functions, variables, logical operators (NOT, OR, AND), quantifiers, parentheses, and some types of equality symbol. There is also a set of rules for recognition of terms and well-formed formulas (wffs).

There are four axioms for quantification:

1) PRED-1: $(\forall x\, Z(x)) \rightarrow Z(t)$
2) PRED-2: $Z(t) \rightarrow (\exists x\, Z(x))$
3) PRED-3: $(\forall x\, (W \rightarrow Z(x))) \rightarrow (W \rightarrow \forall x\, Z(x))$
4) PRED-4: $(\forall x\, (Z(x) \rightarrow W)) \rightarrow (\exists x\, Z(x) \rightarrow W)$

An important theorem for first-order logic is the outcome from Herbrand's work (known as Herbrand's theorem). It states that in predicate logic without equality, a formula *A* in prenex form (all quantifiers at the front) is provable if and only if a sequent *S* comprising substitution instances of the quantifier-free subformula of *A* is propositionally derivable, and *A* can be obtained from *S* by structural rules and quantifier rules only. In other words, it states that the formula

Manuscript created December 29, 2006. Author is Ph. D. student of Department of Information Systems at The Poznan University of Economics, http://www.kie.ae.poznan.pl, email: radekh@teycom.pl.

is provable, if, and only if we can rewrite it without quantifier substituting values and obtain provable formula. For example:

$\forall x\, Z(x) = Z(0) \land Z(1)$

$\exists x\, Z(x) = Z(0) \lor Z(1)$

### C. Boolean Algebra

Boolean algebra (also called Boolean lattice) is an algebraic structure containing objects and operations upon them and set of axioms (see Section D). It consists of one unary operation $\neg$ (not) and two binary operations $\land$ (and), $\lor$ (or) also with two distinct elements 0 (constant representing false), 1 (constant representing true). Language of Boolean algebra considered as language for first order logic also contains symbols: = (equality), $\Rightarrow$ (implication), parentheses and quantifiers, $\forall$ (universal), and $\exists$ (existential).

Boolean algebra has the essentials of logic properties as well as all set operations (union, intersection, complement).

### D. Axioms of Boolean Algebra

Given below is a complete list of Boolean algebra axioms. This set is not a minimal set of axioms (especially staring from Ax13)) – some axioms can be derived from others, but it does not change the reasoning used in this article (the list is larger only for clearness and ensuring that it is complete):

Ax1)  a = b can be written as $(a \land b) \lor (\neg a \land \neg b)$
Ax2)  $a \Rightarrow b = \neg a \lor (a \land b)$
Ax3)  $a \lor (b \lor c) = a \lor b \lor c = (a \lor b) \lor c$
Ax4)  $a \land (b \land c) = a \land b \land c = (a \land b) \land c$
Ax5)  $a \lor b = b \lor a$
Ax6)  $a \land b = b \land a$
Ax7)  $a \lor (a \land b) = a$
Ax8)  $a \land (a \lor b) = a$
Ax9)  $a \lor (b \land c) = (a \lor b) \land (a \lor c)$
Ax10) $a \land (b \lor c) = (a \land b) \lor (a \land c)$
Ax11) $a \lor \neg a = 1$
Ax12) $a \land \neg a = 0$
Ax13) $a \lor a = a$
Ax14) $a \land a = a$
Ax15) $a \lor 0 = a$
Ax16) $a \land 1 = a$
Ax17) $a \lor 1 = 1$
Ax18) $a \land 0 = 0$
Ax19) $\neg 0 = 1$
Ax20) $\neg 1 = 0$
Ax21) $\neg(a \lor b) = \neg a \land \neg b$
Ax22) $\neg(a \land b) = \neg a \lor \neg b$
Ax23) $\neg\neg a = a$

Using universal generalization, one may add in every axiom definition, the universal quantifier stating "for all x *axiom body*" ().

### E. Computational Tree

Corollary of the considerations stated earlier establishes that every formula in Boolean algebra is decidable. It is said to be proved (or called "tautology") if there exists a transformation path from a set of axioms to a sentence that we are trying to prove. Until the authors are discussing FOPC, one may say that every sentence is provable, if, and only if we can start with axiom and repeatedly apply "modus ponens" or "universal generalization" and obtain this sentence [2].

One may then consider every possible (provable) sentence to be deducible from axioms, which may be presented as the graph shown in Fig. 1.

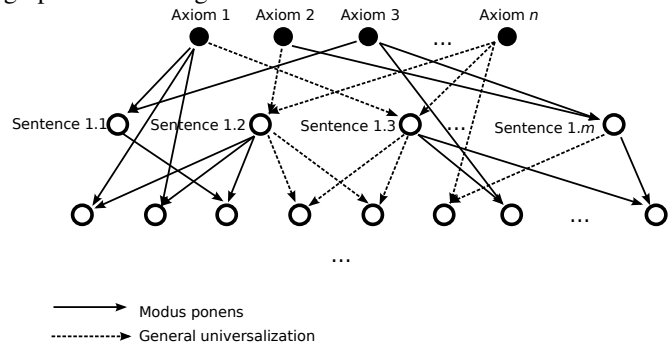

Figure 1 Example of deducible tree

Axioms may also be the result of computations, especially when they are not independent (e.g., ZF axioms) or when computations fall in a cycle. Usually, during computation one would skip deduction to already proven sentences because it does not introduce any new information, so deductions to axioms would have been omitted.

### F. Inference Rules and Deduction

Modus ponens is an inference rule using the reasoning: if a and $a \rightarrow b$ are both proved, then b is also proved.

Universal generalization is an inference rule using the reasoning: if $P(a)$ is proved, and a is a free variable, then $\forall a\, P(a)$ is also proved.

Deduction theorem (in fact deduction meta-theorem) states that if formula F can be deduced from E, then the implication $P \rightarrow Q$ can be directly shown to be deducible from the empty set. Using symbol "$\vdash$" for deducible, one may write: if $P \vdash Q$ then $\vdash P \rightarrow Q$. One may generalize it to a finite sequence of assumption formulas $P_1, P_2, P_3, \ldots, P_n \vdash Q$: $P_1, P_2, P_3 \ldots P_n\text{-}1 \vdash P_n \rightarrow Q$ and repeat it until we obtain the empty set on the left-hand side: $\vdash (P_1 \rightarrow (\ldots(P_n\text{-}1 \rightarrow (P_n \rightarrow Q))\ldots)$.

Deduction follows three kinds of steps: setting up a set of assumptions (hypothesis), reiteration - calling hypothesis made previously to make it recent, and deduction, which is removing recent hypothesis. If one wants to convert proof done using deduction meta-theorem to axiomatic proof, then usually the following axioms would have been involved:

1) $P \rightarrow (Q \rightarrow P)$
2) $(P \rightarrow (Q \rightarrow R)) \rightarrow ((P \rightarrow Q) \rightarrow (P \rightarrow R))$
3) Modus ponens: $(P \land (P \rightarrow Q)) \rightarrow Q$

### G. Corollaries

Theorem 1—If formula expressible in FOPC language is deducible, then every possible transformation of this formula obtained by usage inference rules and axioms is also deducible, and can be expressed in the same language.

Theorem 2—If every transformation of formula is expressible in FOPC, then the optimal for certain resource for chosen computational model is also expressible in the same language.

Proofs of these theorems are provided in the appendix (VI.A and VI.B).

One needs to focus on Theorem 1 and have good understanding of significance of Gödels work. Let there be considered some formula φ which is intended to be proven or disproved. Assuming that there exists some deterministic transformation $T_1$ which transform formula φ to $φ_1$: $φ_1=T_1(φ)$, there can also exist another transformation $T_2$ taking $φ_1$ as input and returning $φ_2$ as output: $φ_2=T_2(φ_1)=T_2(T_1(φ))$. Continuing this idea of transformations one reaches $φ_{TRUE}$ or $φ_{FALSE}$ formula allowing to prove or disprove φ.

Theorem 1 is in fact summary of Gödels Theorem [8], stating that φ, $φ_1$, $φ_2$… $φ_x$ can all be expressed in FOPC language.

Above in fact causes statement of Theorem 2: if every possible transformation / reformulation of formula is expressible in FOPC language, then (by power of every quantification) also optimal transformation is expressible in FOPC language. It does not matter what is nature of this transformations. If they are deterministic (for certain input always returns same output in finite number of steps). Optimal way to solve the problem (decide on formula) can be then written as: $φ_{TRUE/FALSE}=T_x(T_{x-1}(T_{x-1}(…(T_2(T_1(φ)))…))$

### III. cSAT LOWER BOUND

#### A. Problem Definition

In this work, the authors consider a problem called "Count of Satisfaction of Boolean Expression for formula φ." This problem is almost the same as the classical SAT problem, but instead of the question "Is there an assignment to variables such that formula φ is satisfied?", they ask the question "Are there at least $L$ assignments such that formula φ is satisfied?".

$L$ in problem instance is written *unary*, and the remaining part of the instance is exactly the same as in SAT problem (the authors assume that it is in conjunctive normal form (CNF)).

It is easy to show that the problem is in NP – Guess & Check algorithm, for NDTM requires $O(L*n)$ steps to check (certificate size is $L*v$ where $v$ is the number of Boolean variables used).

It is also easy to show that the problem is NP-complete. One can show it using reduction from SAT problem and ask the question "Is there at least $L=1$ assignment such that formula φ is satisfied?"

#### B. Measurable Predicate

Problem question is easy to understand by a human, but it certainly extends to FOPC language defined in Section II. To express it in a defined language, one needs to define predicate "μ" – measure. This predicate is a representation of sigma-additive (countably additive) measurable function known as "set cardinality." Definition of this predicate requires one constant variable $n$ – number of different Boolean variables used. Predicate "μ" will measure number of assignments satisfying formula φ.

1: $μ(∅)$ :- 0
2: $μ(TRUE)$ :- $2^n$
3: $μ(FALSE)$ :- 0
4: $μ(¬φ1)$ :- $2^n-μ(φ_1)$
5: $μ(a1)$ :- $2^{n-1}$
6: $μ(a1∧a2…∧a_k)$
  $∃ a_i, a_j: i≠j ∧ a_i=¬a_j$ :- $μ(FALSE)$
  $∃ a_i, a_j: i≠j ∧ a_i=a_j$ :- $μ(a1∧a2…a_{j-1}∧a_{j+1}…∧a_k)$
    :- $2^{n-k}$
$μ(φ_1∨φ_2)$ :- $μ(φ_1)+μ(φ_2)−μ(φ_1∧φ_2)$

One may think of adding some more conditions to this predicate, but the list given earlier is sufficient to calculate measure for every formula for a defined language (growth of number of axioms and definitions is discussed in Section H). It is also compliant with sigma-measurable function definition.

One may also observe that usage of measure leads to exponential number of calculations required for CNF. This is a consequence of sigma-additivity property: for any sets $a$ and $b$: $μ(a∪b)=μ(a)+μ(b)−μ(a∩b)$. If one considers $m$ sets, then this function transforms to:

$$μ\left(\bigcup_{i=1}^{m} a_i\right) = \sum_{S \in P(\{a_1..a_m\})} \left((-1)^{|S|+1} \cdot μ\left(\bigcap_{a_k \in S} a_k\right)\right), \text{ where}$$

$P(\{a1..a_m\})$ is power set over $m$ sets, which means that it has $2^m$ objects in it. If one is able to calculate the measure of a set or intersection of sets, then the calculation of union of $m$ sets requires $Ω(2^m)$ intersections to be measured.

Problem question using predicate "μ" is then: "$μ(φ)≥L$?".

Direct calculation may not be the only possible way for solving problems, and the authors now analyze the definition of lower bound, deterministic and nondeterministic computation models.

#### C. Lower Bound Definition

Lower bound in Big-O notation is denoted as $Ω(g(n))$, and for its use in this article, one may assume that it is used to express *problem lower bound*. Interpretation of lower bound is "minimum value of function in the worst case," and is defined as $f(n) \in Ω(g(n)) \Leftrightarrow \liminf_{n \to \infty} \left|\frac{f(n)}{g(n)}\right| > 0$.

In most of the complexity considerations, two types of resources are used in the expressions of problem lower bounds or algorithm upper bounds. These resources are time (number of steps required) and space (number of symbols/tape cells required).

Theorem 3—Time complexity of problem/algorithm is always greater than or equal to space complexity. This theorem is proved in Section VI.C.

Theorem 4—Minimal number of symbols required for unambiguous description of object is $Ω(\log(N))$, where $N$ represents the number of possible objects to be stored. In other words, this means that if one has $N$ different objects that may occur in computations at a certain step and would want to store information on which one occurred, then $Ω(\log(N))$ symbols are required. This theorem is proved in Section VI.D.

In this work, the authors mainly consider time complexity using observation from Theorems 3 and 4.

Theorem 5—Lower bound calculated to express a specific resource (time or space) usage for deciding formula expressed in FOPC for a chosen computational model is equal to the minimal usage of this resource for the best possible transformation of formula in this language. This theorem is proved in Section VI.E.

Theorem 5 is consequence of Theorem 2. If one had set of deterministic transformations expressing optimal way to solve the problem:

$\varphi_{TRUE/FALSE}=T_x(T_{x-1}(T_{x-1}(\ldots(T_2(T_1(\varphi)))\ldots)))$ then by power of definitions it can be shown that lower bound for problem solution is exactly equal to time required by this optimal solution. This is consequence of lower bound definition – it is asymptotically minimal amount of resource required to solve the problem. Repeating most important observations till this point:

a) formula $\varphi$ can be expressed in FOPC language (from Gödels Theorem [8])
b) any possible transformation of formula can be expressed in FOPC language $\varphi_1=T_1(\varphi)$ (Theorem 1)
c) if every deterministic transformation can be expressed in FOPC language then also optimal deterministic transformation can be expressed in FOPC language (Theorem 2)
   $\varphi_{TRUE/FALSE}=T_x(T_{x-1}(T_{x-1}(\ldots(T_2(T_1(\varphi)))\ldots)))$
d) resource cost of optimal transformation of formula is equal to deterministic lower bound of the problem

Roughly speaking, lower bound should be considered as the minimal amount of resource used for computation for the worst case. In case of time, it is the minimal number of operations to perform. It is even intuitive to see that if one could express calculation in some "steps," then lower bound is equivalent to minimal number of "steps" required in the worst case.

### D. Nondeterministic Calculation Model

Nondeterministic calculation model may be considered as the "luckiest possible guesser." Such an approach expresses that the role of NDTM to answer a problem question is to guess the certificate and check it. If the check can be performed in $O(n^c)$ for some constant $c$, then one considers the problem as part of NP complexity class.

One has to remember that DTM is a "special case" of NDTM where from every machine state, only one possibility to choose the next state exists, regardless of the symbol in the cell where the tape read/write head is positioned. This means that every problem solvable by DTM in $O(n^c)$ steps is solvable also on NDTM in at most same number of steps (or may be less).

A good example expressing the differences between DTM and NDTM is the 2SAT problem (classic satisfaction of Boolean expression in CNF problem, but where in each clause there are at most two literals). This problem is solvable by DTM in $O(n3)$ steps, but NDTM may guess the correct assignment and verify it in $O(n)$.

In terms of first-order logic and Herbrand's theorem, one can see that NDTM is a verifier of Herbrand's subformulas. When the formula is expressed using existential quantifier: $\exists$ <a: assignment> $F(a)$, then, according to Herbrand's theorem, it is equivalent to: $F(a1) \vee F(a2)\ldots \vee F(a_k)$. NDTM is able to check each of $F(a_i)$ simultaneously, even if the number of possible assignments is exponential, excepting when at least one of the computation paths led to an accepting state.

One can see that for a nondeterministic calculation model problem, the number of steps of lower bound is equal to the minimal number of steps required to check the certificate.

For example, for 2SAT problems one can have different approaches.

| App | Number of possible Herbrand's subformulas | Minimal number of steps to check each subformula | Total calculation cost |
|---|---|---|---|
| 1 | $2^n$ | $N$ | $n$ |
| 2 | $2^p$ (guessing only $p$ variables) | $n*p^3$ | $N*p^3$ |
| … | | | |
| K | 1 (without splitting) | $n^3$ | $n^3$ |

Table 1 Different approaches for nondeterministic calculation

Table 1 presents different approaches differing mainly in the number of "guesses." Calculation of problem lower bound for nondeterministic model of calculation returns the minimal number of steps required to check subformula.

The last row presents the approach where the problem is not split, so it is calculated as in the deterministic model of calculations.

### E. Deterministic Calculation Model

As mentioned in Section D, deterministic model of calculations follows a single computation path. It is obvious that despite direct calculations, DTM can also perform Guess & Check algorithm (simulating NDTM). This time, the authors do not assume that DTM is the "luckiest possible guesser" and for lower bound complexity calculation of this approach, they have to assume that DTM is the "worst possible guesser." This is also a consequence of the slight change in computation goal - NDTM has to "accept" when there is computational path leading to accepting state, while DTM has to "decide" on input, which means that the answer "NO" can be produced only when there is no possible way of reaching the accepting state (NDTM can be defined without rejecting state).

Additionally, DTM requires an iterator (space on tape where number of current "guess" can be stored), which according to Theorem 4 requires $\Omega(\log(H))$, where $H$ represents the possible number of "guesses."

Table 2 shows what time complexity would look like.

| App | Number of possible Herbrand's subformulas | Minimal number of steps to check each subformula | Total calculation cost |
|---|---|---|---|
| 1 | $2^n$ | $N$ | $2^n*n+\log(2^n)$ |
| 2 | $2^p$ | $n*p^3$ | $2^n*n*p^3+\log(2^p)$ |
| … | | | |
| K | 1 | $n^3$ | $n^3$ |

Table 2 Different approaches for deterministic calculation

In this table, the last row represents the minimal possible number of steps to calculate result. It is easy to show that for DTM, this row also presents deterministic problem lower bound because if any of the "guessing" approaches had been better, then it would have been used to present minimal deterministic calculation cost (we assume that values in the table are "best possible" not "best known" - see Theorem 5).

### F. cSAT Nondeterministic Algorithm Upper Bound

Upper bound for algorithm solving cSAT problem is polynomial. It is a consequence of Herbrand's theorem and ability of NDTM to:

generate all possible subformulas in O($n^c$)

verify each of them in O($n^c$)

NDTM algorithm can be described using the following steps:
1) Guess sets of measure $L$ consisting of assignments of variables (time O($L*v$))
2) Verify guessed set (time O($L*v$))

This procedure leads to accepting the state (if at least one computation path is accepting) in at most O($L*v$) steps and because instance size $n \in \Omega(L+v)$, the solution is provided in O($n2$).

### G. cSAT Deterministic Lower Bound

Now, using the observations described in the earlier sections, the authors calculate deterministic lower bound of cSAT problem (it is known that its nondeterministic upper bound is O($L*v$)).

First, one needs to write the problem in the FOPC language. One uses the predicate µ: µ(φ)≥$L$.

This problem may be considered to be harder than the classic SAT problem. If one tries to guess all possible subsets, then we would have $\Omega(2^{2^v})$ subsets, so according to Theorem 4, it would require $\Omega(2^v)$ symbols to store information about the considered subset, which, according to Theorem 3, leads to the conclusion that such a calculation requires at least $\Omega(2^v)$ steps. "Guessing" only subsets of size $L$ leads to $\Omega(2^v)$ different subsets, so that it can be calculated by NDTM in polynomial time ($\Omega(v*L)$ steps), but requires $\Omega(2^v*v*L)$ steps to calculate on DTM.

In fact, following the assumption that DTM is the worst possible guesser, one may see that the number of hypotheses ("guesses") used during computation can lead to an exponential usage of time if "depth" of hypothesis path is longer than O(log($n$)) or any of the hypotheses has more than polynomial number of possible values. For example, if one states hypothesis A with possibilities, it is true or false (constant number of possibilities) and it is followed by hypothesis B (true or false), etc.; we need O($n$) hypotheses before we can decide on formula, then in the worst case we require $\Omega(2^n)$ steps to give the answer "NO."

Leaving then all Guess & Check approaches, the authors try to determine the minimal possible number of steps for DTM to decide on problem input. According to Theorem 5 and considerations from the earlier sections, the authors conclude that the shortest possible path consists of steps transforming input formula to axioms of theory. If one can show that *every* transformation requires exponential number of steps or usage of object using exponential number of symbols to store, then it will be direct proof that lower bound of cSAT problem is over-polynomial.

When will one be able to observe exponential growth of minimum number of required steps? If after using an axiom or predicate, one will obtain a formula of multiplicative length by a factor greater than 1. For example, if for formula of size $n1$ (considered to be in CNF), the authors use Ax9) for one parenthesis, they obtain new formula in the format $v1\wedge(n2,1)\vee v2\wedge(n2,2)\vee…\vee v_m\wedge(n2,_m)$. In each of the $m$ parts, one can use a variable from the beginning to remove all its negations from body, so |$n2,_*$|<|$n1$|–$m$, but for very large $n1$, these parts of formula will still require further transformations, which if done only with Ax9) would lead to exponential growth. Concluding this paragraph, one may say that if transformation reduces size of formula substring by O($n^c$) and multiplies this shorter string in formula making string grow to $n2$, where $n2 \in \Omega(n1*c)$, then this path leads to exponential growth of formula and thus its lower bound is $\Omega(2^n)$.

In Table 3, the authors present the effect obtained by usage of every possible transformation, but before this the authors define polynomial *purifying* function for formula. This function will use axioms Ax7), Ax8), Ax11), Ax12), Ax13), Ax14), Ax15), Ax16), Ax17), Ax18), Ax19), Ax20), Ax23), and two observations:

µ(φ$_1$)=µ(φ$_1 \wedge$(T$_{RUE}$))=µ(φ$_1 \wedge$($v1 \vee v2 \vee … \vee$T$_{RUE}$));

µ($v1 \wedge$φ1)=µ(φ$_2$) where φ$_2$ is obtained by replacing every occurrence of $v1$ in φ1 with constant TRUE.

Roughly speaking, this function looks for variables that can be cleared out from formula and prepare it for the next step of calculation. The authors assume that at every step of calculation, formula is in a form not allowing the use of any of the above axioms or rules. It is also important to remember that the number of transformation rules does not matter - refer to Section H.

| Transformation used | Length of string used | Result string length | Lower bound for path | Remarks for "worst case" |
|---|---|---|---|---|
| Ax1) | | | | These axioms cannot be used since input never contains these symbols |
| Ax2) | | | | |
| Ax3) | $n_1$ | $n_1$ | $\Omega$(cSAT) | These axioms do not change formula length |
| Ax4) | $n_1$ | $n_1$ | $\Omega$(cSAT) | |

| Transformation used | Length of string used | Result string length | Lower bound for path | Remarks for "worst case" |
|---|---|---|---|---|
| Ax5) | $n_1$ | $n_1$ | $\Omega(cSAT)$ | |
| Ax6) | $n_1$ | $n_1$ | $\Omega(cSAT)$ | |
| Ax7) | | | | These axioms cannot be used because formula is transformed by *purifying* function |
| Ax8) | | | | |
| Ax9) | $m_1+m_2+n_r$ | $2*m_1*m_2+n_r-2$ | $\Omega\left(\frac{n_1}{p_1}*\overline{m_1}^{p_1}\right)$ | Used on two parenthesis replaces them with string of size $2*m_1*m_2$, after using these axioms *purifying* function will reduce size but in the worst case only by 2 symbols |
| Ax10) | $m_1+m_2+n_r$ | $2*m_1*m_2+n_r-2$ | $\Omega\left(\frac{n_1}{p_1}*\overline{m_1}^{p_1}\right)$ | |
| Ax11), Ax12), Ax13), Ax14), Ax15), Ax16), Ax17), Ax18), Ax19), Ax20) | | | | These axioms cannot be used because formula is transformed by *purifying* function |
| Ax21) | $n_1$ | $n_1$ | $\Omega(cSAT)$ | These axioms do not change formula length |
| Ax22) | $n_1$ | $n_1$ | $\Omega(cSAT)$ | |
| Ax23) | | | | As Ax20) above |
| μ1 | | | | These rules cannot be used because formula is transformed by *purifying* function |
| μ2 | | | | |
| μ3 | | | | |
| μ4 | $n_1$ | $n_1$ | $\Omega(cSAT)$ | Can be used with Ax21), Ax22) or Ax23), but does not change length of formula |
| μ5 | | | | As μ3 above |
| μ6 | $m_1+m_2$ | $m_1*m_2-2$ | $\Omega\left(\frac{n_1}{p_1}*\overline{m_1}^{p_1}\right)$ | Treating formula as consisting of two parts |

Table 3 Lower bounds for every possible transformation of cSAT formula

The variables used in Table 3 are the following:
$n1$ – Length of the formula
$n_r$ – Length of the remaining part of the formula
$m1$ – Length of first part/parenthesis of the formula
$m2$ – Length of second part/parenthesis of the formula
$p1$ – Number of parentheses in the formula
$\overline{m_1}$ – The authors consider asymptotic behavior of the function, so one may use kind of "mean" $m_1$ – representing $\Omega(m_1)$.

It is clear that in the worst case, Ax9), Ax10), or μ6 have to be used several times before *purifying* function would make significant reduction of length. Lower bound is considered as the *minimal worst case,* so from this table it is clear that in the worst case it is $\Omega(m^p)=\Omega(2^p*\log(m))$ and because $p$ and $m$ are both O($n$), the whole lower bound is $\Omega(2^n)$.

### H. More Conditions and More Axioms

The above considerations prove clearly that deterministic lower bound for cSAT problem considered with FOPC language defined is exponential. But one needs to answer one additional question – Is it the result of too poor FOPC axioms set definition? Or are too few predicates defined?

In [1] Baker–Gill–Solovay theorem, authors have shown that problem "Is P equal to NP?" can be relativized using oracles. Oracle is a machine (black box) that gives answers to certain type of problems in one step. One can then imagine that there are a very large number of oracles which can solve certain types of instances. DTM task is to pick up one of them (or use them sequentially because if the number of oracles is an attribute of the machine, then even if we have used millions of them, the complexity in terms of relation to instance size is O(1)).

The authors presume then, that for cSAT problem, there exists some deterministic algorithm calculating answer in O($n^c$) steps. Following lower bound calculation, one knows that this algorithm calculates a result requiring $\Omega(2^n)$ transformations. Reminding optimal transformation as described above: $\varphi_{TRUE/FALSE}=T_x(T_{x-1}(\ldots(T_2(T_1(\varphi)))\ldots))$ and $x\in\Omega(2^n)$.

The presumption made here can be presented as existence of some transformation $T_A\equiv T_{x-k}(T_{x-k-1}(\ldots(T_{x-k-m}(\quad))\ldots))$, where $m$ is exponential ($T_A$ is equivalent to exponential number of transformations in FOPC language, on optimal transformation path). The authors also assume that $T_A$ is deterministic, as deterministic lower bound is discussed in this section, and computable in polynomial number of steps.

Now, one need to look on transformation path as on decision process, where at each step there is a decision to be made (decision which transformation is to be used). Each decision takes $\Omega(1)$ space to be stored. If $m$ was dynamic and asymptotically equal to $2^n$, and also computable in polynomial number of steps then this would be equal to O($2^n$) decisions in O($n^c$) time what contradicts Theorem 3.

Considering constant $m$ (invented by algorithm designer) one may ask what is common in a large number of Turing machines (in the sense of defined algorithms), large number of axioms, large number of predicates, large number of oracles, or large $m$ in above transformation $T_A$? Their number is always a *constant*, even if very large. If then anyone defines multiple TMs, adds multiple axioms, predicates, defines large number of oracles, or finds one transformation $T_A$ equivalent to exponential number of other transformation, then in fact after defining them, one may have constant number of machines, axioms, predicates, transformations, and oracles.

The authors then assume that there exists a machine denoted by LDTM in which implements are equivalent to large number of TMs, large number of axioms, predicates, implements $T_A$ and are connected to multitude of oracles. Such a defined machine is (by power of assumption) capable of answering cSAT questions for a finite number of differing input types (number of types is a consequence of maximal input size).

In other words, the authors assume that there exists a machine LDTM able to answer cSAT questions for instance size less than or equal to $n_l$. They may consider having $O(c^{n_l})$ different input types, and each type is covered at least by one combination of axioms, predicates, or oracles allowing LDTM to give answer in $O(n)$ steps. One may assume that there are $g_l$ such combinations. Denoting $|g_i(n_l)|$ number of instances solved by $i$th combination of axioms and oracles for instance $n_l$ symbols long, the authors have assumed that: $\sum_{i=1}^{l}|g_i(n_l)| \geq c^{n_l}$, where $\forall_{i=1}^{l} |g_i(n_l)| \lessdot c^{n_l}$.

Now the authors determine the ability of this machine to answer cSAT question where $n=n_l*y$. Number of combinations of axioms and oracles remain constant ($g_l$), but they assume that each combination covers more instances (considered to be "same type") $g_l(n)=g_l(n_l*y) \leq g_l(n_l)^y$. The number of possible types grows from $O(c^{n_l})$ to $O(c^{n_l*y}) = O((c^{n_l})^y)$.

Calculating instances covered by $g_l$ definitions, we have: $\sum_{i=1}^{l}|g_i(n_l*y)| \leq \sum_{i=1}^{l}|g_i(n_l)|^y$.

If one proves that for $y$ growing to infinity $\sum_{i=1}^{l}|g_i(n_l)|^y \lessdot c^{n_l}|^y$, it will be proof that not all instances of size $O(n_l*y)$ are solvable using LDTM definitions, so these large instances will require calculations using deterministic lower bound discussed in Section G. Proof is presented in VI.F, so corollary about impossibility to answer cSAT problems in polynomial time by LDTM holds.

## IV. COROLLARIES

To summarize this article, the authors repeat the deduction path:
a) formula φ can be expressed in FOPC language (from Gödels Theorem [8])
b) any possible transformation of formula can be expressed in FOPC language $\varphi_1=T_1(\varphi)$ (Theorem 1)
c) if every deterministic transformation can be expressed in FOPC language then also optimal deterministic transformation can be expressed in FOPC language (Theorem 2)
$\varphi_{TRUE/FALSE}=T_x(T_{x-1}(T_{x-1}(\ldots(T_2(T_1(\varphi)))\ldots)))$
d) resource cost of optimal transformation of formula is equal to deterministic lower bound of the problem
e) $T_A$ equivalent to exponential number of transformations computable in polynomial time contradicts Theorem 3
f) large number of defined constant set of transformations, oracles, algorithms, machines ect. cannot cover all possible inputs for growing instance size (Theorem 6)
g) optimal solution of problem requires $\Omega(2^n)$ transformations
h) deterministic lower bound for cSAT problem is then $\Omega(2^n)$, then cSAT$\notin$P (Theorem 5)
i) NDTM solves cSAT in polynomial time, so cSAT$\in$NP
j) this means that P$\neq$NP

If above considerations are correct then checking problem known to be in P has to show that it is in P using the same reasoning. Such check for 2SAT problem if presented in Section VI.G - lower bound for this problem is $\Omega(n^c)$.

In [1], there was presented an oracle A for which $P^A=NP^A$. Proof presented in Section VI.F and problem lower bound lead to corollary, and if A is able to solve cSAT in polynomial time, then A has to be nondeterministic - it has to consist of infinite number of objects: deterministic oracles, algorithms, DTMs, axioms, rules etc. (the authors also consider NDTM as an infinite set of DTM duplicates - each for one computational branch).

This work discusses problem P=NP, as described in [5]. It may be said to *relativize* (see [1]) to deterministic model of computation showing that deterministic calculation model made up of finite number of machines (algorithms), oracles, axioms, or predicates is incapable of solving the considered problem when its instance grows to infinity.

On the other hand, one may conclude that if restrictions on maximum input length problem are set, then the problem can be proved to be in P using a large number of machines, axioms, algorithms, predicates, or oracles.

For deterministic model of computation, one knows then that P$\neq$NP.

Using Theorem 13 from [12], the authors also know that NP-complete$\neq$(NP-P). In this theorem, the authors have proved that: if P$\neq$NP and U is some NP-complete language then U=A$\cup$B where neither A nor B language is NP-complete (at least one of them is also not equal to P: A$\neq$P $\vee$ B$\neq$P). Complexity classes can be put in a picture (Fig. 2):

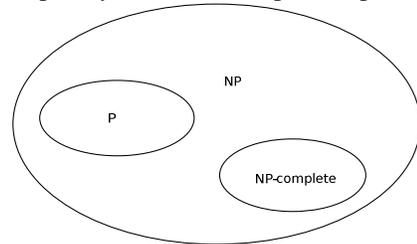

Figure 2 Relation between P, NP, and NP-complete classes

## VI. APPENDIX

### A. Proof 1 - Proof of Theorem 1

Theorem 1 - If formula expressible in FOPC language is deducible, then every possible transformation of this formula obtained by usage inference rules and axioms is also deducible and can be expressed in the same language.

This theorem is a direct consequence of FOPC definitions. If $\varphi$ is deducible, then:

- $\varphi \wedge$ axiom
- $\varphi \rightarrow$ axiom
- $\forall x\, \varphi$
- …

are also deducible.

### B. Proof 2 - Proof of Theorem 2

Theorem 2 - If every transformation of formula is expressible in FOPC, then the optimal for certain resource for chosen computational model is also expressible in the same language.

This theorem is a consequence of Theorem 1 and FOPC definitions. If the goal of calculation is to decide on formula based on theory axioms, then it is required to obtain formula as a consequence of axioms (with empty left-hand side):

$\vdash (P_1 \rightarrow (\ldots(P_{n-1} \rightarrow (P_n \rightarrow Q))\ldots))$.

The authors said that *every* possible transformation of formula is expressible in FOPC and this directly means that the optimal in the aspect of a certain resource (time or space) path is also expressible in FOPC.

### C. Proof 3 - Proof of Theorem 3

Theorem 3 - Time complexity of problem/algorithm is always greater than or equal to space complexity.

This theorem is a consequence of Turing machine definition, which states that in one step, a machine can read or write one (or in general constant) number of symbols. If then $f(n)$ symbols were written, then machine had used at least $f(n)$ steps to write them.

### D. Proof 4 - Proof of Theorem 4

Theorem 4 - Minimal number of symbols required for unambiguous description of object is $\Omega(\log(N))$, where $N$ represents the number of possible objects to be stored.

In this section, the function log is considered to have $\Sigma$ in root, where $\Sigma$ represents the number of the symbols in the alphabet: $\log(\Sigma)=1$.

The authors prove the theorem using contradiction. Suppose that one knows "compression" algorithm allowing to write each of $N$ symbols using $\log(N)-f(N)$ symbols, where $f(N)$ is a function such that: $\forall N: 0<f(N)<\log(N)$.

On $\log(N)-f(N)$, one can write at most $\Sigma^{\log(N)-f(N)}$ different strings.

$$\Sigma^{\log(N)-f(N)} = \frac{\Sigma^{\log(N)}}{\Sigma^{f(N)}} = \frac{\Sigma^{\log(N)}}{\Sigma^{f(N)}} = \frac{N^{\log(\Sigma)}}{\Sigma^{f(N)}} = \frac{N}{\Sigma^{f(N)}}$$

Now the authors check whether this number is greater than the number of objects to be identified by checking limens:

$$\lim_{N \to \infty}\left(\frac{N}{\Sigma^{f(N)}} - N\right)$$

It is easy to see that if $f(N)=0$, then limens is equal to zero (which means that exactly $N$ different objects can be described using a string of desired length), but when $f(N)>1$ (it is the smallest value making difference in the number of symbols used), it is negative which means that less than $N$ objects can be represented using a string of this length.

### E. Proof 5 - Proof of Theorem 5

Theorem 5 - Lower bound calculated to express specific resource (time or space) usage for deciding formula expressed in FOPC for a chosen computational model is equal to the minimal usage of this resource for best possible transformation of formula in this language.

Proof of this theorem is in fact a direct corollary of Theorems 1 and 2. If any transformation of formula is expressible in FOPC language, then the optimal in terms of chosen resource is also expressible in FOPC language and when the authors calculate lower bound for this resource for the whole transformation path (from input string to decidable form (to axioms)), then they obtain the value of lower bound for the considered problem.

### F. Proof 6 - Proof of Not Covering by Constant Set of Definitions All Possible Large Instances by LDTM

Assumptions: $\forall_{i=1}^{l} |g_i(n_l)| < |c^{n_i}|$ and $\sum_{i=1}^{l} |g_i(n_l)| \geq |c^{n_i}|$. Also $g_i(n)$ function operates on natural numbers and returns natural numbers, so $\forall_{i=1}^{l} |g_i(n_l)| \leq |c^{n_i} - 1|$.

One want to solve $\sum_{i=1}^{l} |g_i(n_l)|^y < |c^{n_i}|^y$ for $y$ growing to infinity.

First one may observe that: $\sum_{i=1}^{l} |g_i(n_l)|^y \leq \sum_{i=1}^{l} |c^{n_i} - 1|^y$ if one can solve inequality $\sum_{i=1}^{l} |c^{n_i} - 1|^y < |c^{n_i}|^y$, then it will be equivalent to prove that the proof is correct.

The new equality presented by the authors is free from $i$ variables, so it can be rewritten as: $l * |c^{n_i} - 1|^y < |c^{n_i}|^y$.

Now the authors take logarithm on both sides to the base $l$:
$\log_l(l * |c^{n_i} - 1|^y) < \log_l(|c^{n_i}|^y)$
$\log_l(l) + \log_l(|c^{n_i} - 1|^y) < \log_l(|c^{n_i}|^y)$
$1 + y * \log_l(|c^{n_i} - 1|) < y * \log_l(|c^{n_i}|)$
$\frac{1}{y} + \log_l(|c^{n_i} - 1|) < \log_l(|c^{n_i}|)$

At this point, it is obvious that this inequality holds - $1/y$ when $y$ grows to infinity may be omitted and one has inequality of two logarithms where this one on the left-hand side is the logarithm of lower value.

More formally, one may calculate limens:
$\lim_{y \to \infty} \left( \frac{1}{y} + \log_l(|c^{n_i} - 1|) - \log_l(|c^{n_i}|) \right)$
$= \lim_{y \to \infty} \left( \log_l \left( \frac{|c^{n_i} - 1|}{|c^{n_i}|} \right) \right) = \lim_{y \to \infty} \left( \log_l \left( 1 - \frac{1}{|c^{n_i}|} \right) \right) = -x$

Proof is then correct.

### G. Proof 7 - Lower Bound for 2SAT Problem

2SAT problem is a special case of cSAT problem. Its special factors are:
- $L = 1$ (problem question is "$\mu(\varphi) \geq 1$" or "$\mu(\varphi) > 0$")

$m_1 = m_2 = \ldots = m_p = 2$

The authors assume that the input string is in CNF and *purifying* function (defined in Section III.G) cannot be applied.

They use Ax9) on parenthesis to select next parenthesis such that:
- parenthesis has not been used yet
- parenthesis contains negation of variable used in a previous step

Every time the usage of Ax9) will be followed by *purifying* function usage.

For example:
$(a \vee b)1 \wedge (a \vee \neg c)2 \wedge (c \vee d)3 \wedge (\neg b \vee \neg c)4 \wedge (\neg b \vee \neg a)5 \wedge (\neg d \vee \neg a)6 \wedge (e \vee f)7$ the authors would have used Ax9) for parenthesis: 1 and 6 (last variable is: $\neg d$), then 3 (lv: $c$), then 4 (lv: $\neg b$), then 1 (second time, lv: $a$), then 5 (lv: $\neg b$), then 1 (third time, lv: $a$), then 6 (second time, lv: $\neg d$), then 3 (second time, lv: $c$), and finally 2. Every parenthesis will be used at most on every path from any pair of parentheses.

At every stage, calculation formula will contain at most $p+1$ conjunctions (where $p$ is the number of parenthesis processed) and each conjunction will contain at most every variable once. Every parenthesis will be used at most $p2$ times, which means that at every stage of computation, formula length is $O(p3)$.

This may not be a time optimal solution. According to Theorem 2, optimal transformation path is expressible using axioms and predicates defined for FOPC, but to show that the problem is in P, one does not need to look for optimal transformation path - the authors have shown that there exists at least one transformation path polynomially bounded to instance size, and even if it is not the optimal one, it shows that 2SAT problem is in P.